\documentclass[12pt,leqno]{article}
\newcommand{\EE}{{\rm I\kern-2pt E}}
\newcommand{\RR}{{\rm I\kern-2pt R}}
\newcommand{\DD}{{\rm I\kern-2pt D}}
\newcommand{\PP}{{\rm I\kern-2pt P}}
\newcommand{\NN}{{\rm I\kern-2pt N}}
\newcommand{\dd}{{\rm \kern 3pt I\kern-9pt d}}

\renewcommand{\thesubsection}{\thesection.\arabic{subsection}}

\title{\large ERROR CALCULUS AND PATH SENSITIVITY IN FINANCIAL MODELS}
\author{\sc Nicolas Bouleau}
\date{\it Ecole des Ponts, ParisTech}

\begin{document}
\maketitle

{\small
In the framework of risk management,
for the study of the sensitivity of pricing and hedging in stochastic financial models to changes of parameters
 and to perturbations of the stock prices, we 
propose an error calculus which is an extension of the Malliavin calculus based on Dirichlet forms. Although useful also in physics, this error calculus is well
adapted to stochastic analysis and seems to be the best practicable in finance.
	
This tool is explained here intuitively and with some simple examples.\\

}\\
\setcounter{section}{1}
\begin{center}
{\thesection. INTRODUCTION}
\end{center}

Once a model is chosen to price contingent claims and to hedge a position, the actual
 questions obviously are : what is the exposure to errors on the model and to changes in the market ? This risk assessment is usually
done in terms of sensitivity of the portfolio to variations of numerical financial quantities and  parameters of the model. 
Now the theory of Dirichlet forms allows to extend this sensitivity calculus to perturbations of  functional quantities like stochastic processes.
The method is an extension of Malliavin calculus. The errors are thought to be infinitesimal random quantities with biases, 
 variances and covariances.  

Among the promising clues of application of error calculus to finance, let us mention some directions : to manage the precision of numerical methods used
to implement the stochastic theory ; to obtain new integration by parts formulas to speed up Monte Carlo simulations ; to study how depends the solution
of an ODE or of an SDE on a functional coefficient like the level-type volatility. We sketch these two last directions at the end of the article.

Our presentation starts from the basic ideas and goes until examples of completely tractable computations.

\setcounter{section}{2}
\begin{center}
{\thesection. PRESENTATION OF THE METHOD}
\end{center}
\setcounter{equation}{0}
\setcounter{subsection}{1}
\label{sec:2.2}
Several technics are available to represent errors mathematically and to compute them with formulas. The method of Dirichlet forms is based on some major
historical ideas which are the easiest way to penetrate it.\\

\noindent{\thesubsection. Propagation of errors : error calculus \`{a} la Gauss}\setcounter{subsection}{2}\\

\noindent After his argument showing the importance of the normal law ({\it Theoria
 motus corporum
coelestium} 1809), Gauss
 was interested in the propagation of errors 
({\it Theoria Combinationis} 1821). Given a quantity
 $U=F(V_1, V_2,\ldots)$ function of other erroneous quantities $V_1, V_2,\ldots$ he states the problem
 of computing the quadratic error to fear on $U$ knowing the quadratic errors $\sigma_1^2,\sigma_2^2, \ldots$ on $V_1, V_2,\ldots,$ these
errors being supposed small and independent. His answer is the following formula
\begin{equation}
\sigma_U^2=(\frac{\partial F}{\partial V_1})^2\sigma_1^2 +
(\frac{\partial F}{\partial V_2})^2\sigma_2^2 + \cdots
\end{equation} 
\noindent he gives also the covariance between the error on $F$ and the error of an other function of the $V_i$'s.

Formula (2.1) possesses a property which makes it highly better, in several questions, than other formulas
 used here and there in textbooks during the 19th and 20th centuries. It is a {\it coherence} property. By lack of place we refer
 to Bouleau (2001) \S 1 for the reason of this coherence property.

In the calculus \`{a} la Gauss the errors on $V_1,V_2,\ldots$ are not necessarily supposed to be independent nor constant, they can depend on $V_1,V_2,\ldots$ 
: Let be given a field of symmetric positive matrices $(\sigma_{ij}(v_1,v_2,\ldots))$ on ${\RR}^d$ representing
 the conditional variances 
and covariances of the errors on $V_1,V_2,\ldots$ given the values $v_1,v_2,\ldots$ of $V_1,V_2,\ldots$ then the error on $U=F(V_1,V_2,\ldots)$
is \begin{equation}{\textstyle
\sigma_F^2=\sum_{ij}\frac{\partial F}{\partial V_i}(v_1,v_2,\ldots)\frac{\partial F}{\partial V_j}(v_1,v_2,\ldots)\sigma_{ij}(v_1,v_2,\ldots)
}\end{equation}
which depends solely on $F$ as mapping.\\

\noindent{\thesubsection. First order and second order calculus}\setcounter{subsection}{3}\\

The following remark, although very simple, is important to understand the role of the error calculus \`{a} la Gauss that  will be used in the sequel
in the extended form allowed by Dirichlet forms.

Let us start with a quantity $x$ with a small centred error $\varepsilon Y$, on which acts a non-linear regular function $f$.
Thus we have at the beginning a random variable written $x+\varepsilon Y$, it has no bias (centred at the true value $x$)
and its variance is $\varepsilon^2\sigma_Y^2$.

After having applied the function $f$, using Taylor formula shows that the error is no more centred. The bias has the same
 order of magnitude as the variance.
Then applying  new regular non-linear functions $f_n$ gives  a transport formula which shows how errors propagate : biases and variances keep permanently
the same order of magnitude
$$\begin{array}{rcl}
{\mbox{bias}}_{n+1}&=&{\mbox{bias}}_n f'_{n+1}(x_n)+\frac{1}{2}{\mbox{variance}}_n f''_{n+1}(x_n)+\varepsilon^30(1)\\
{\mbox{variance}}_{n+1}&=&{\mbox{variance}}_n f{'^2}_{\!\!\!\!\!n+1}(x_n)+\varepsilon^30(1)
\end{array}
$$
(it could be easily extended to applications from ${\RR}^p$ to ${\RR}^q$, for the general formulas on the bias and the variance
of the error under regular mappings see Bouleau and Hirsch (1991) chapter I \S 6 corollaries 6.1.3 and 6.1.4).

We see that {\it the calculus on the biases is a second order calculus involving the variance. Instead, the calculus on the variances
 is a first order calculus not involving the biases.} Surprisingly, the calculus on the second order moments of errors is indeed simpler than the calculus on the
first moments. 
Thus, the error calculus on the variances appears to be necessarily the first step in an analysis of errors propagation based on differential methods
and supposing small errors.\\

\noindent{\thesubsection. Extended error calculus using Dirichlet forms}\setcounter{subsection}{4}\\

\label{subsec:2.3}
The error calculus of Gauss has the limitation that it has no mean of extension. If the error on 
$(V_1,V_2,V_3)$ is known it gives the error on any dif\-ferentiable function of $(V_1,V_2,V_3)$
but that's all.

Now, in the usual probabilistic situations where 
a sequence of quantities $X_1,X_2,\ldots,X_n,\ldots$ is given and where the errors are known on the regular functions of a finite number
of them, we would like to deduce the error on a function of an infinite number of the $X_i$'s or at least on some such functions.

It is actually possible to {\it reinforce} this error calculus giving it a powerful extension tool and preserving the coherence property. In addition,
 it will give us the comfortable
 feature to handle Lipschitz functions as well.

For this we come back to the idea that the erroneous quantities are themselves random, as Gauss had supposed for his proof of the 
`law of errors', say defined on $(\Omega, {\cal A}, {\PP})$. The quadratic error on a random variable $X$ is then itself a random variable that
we will denote  by $\Gamma[X]$. Intuitively we still suppose the errors are infinitely small although this doesn't appear in the notation.
It is as if we had an infinitely small unit to measure errors fixed in the whole problem. The extension tool is the following, we 
assume that if $X_n \rightarrow X$ in $L^2(\Omega, {\cal A}, {\PP})$ and if the error $\Gamma[X_m-X_n]$ on $X_m-X_n$ can be made as small as
we want in $L^1(\Omega, {\cal A}, {\PP})$ for $m,n$ large enough, then the error $\Gamma[X_n-X]$ on $X_n-X$ goes to zero in $L^1$.

This can be axiomatized as follows : we call {\it error structure} a probability space equipped with a local Dirichlet form possessing a 
{\it carr\'{e} du champ}.

 \noindent{\sc Definition 2.1.} {\it An error structure is a term
$(\Omega, {\cal A}, {\PP}, {\DD}, \Gamma)$
where $(\Omega, {\cal A}, {\PP})$ is a probability space, satisfying the four properties :

1.) ${\DD}$ is a dense subvectorspace of $L^2(\Omega, {\cal A}, {\PP}).$

2.) $\Gamma$ is a positive symmetric bilinear map from ${\DD}\times{\DD}$ into $L^1({\PP})$ fulfilling the functional calculus of
 class ${\cal C}^1\cap {\mbox{Lip}}$, what means that if $u\in{\DD}^m$ and $v\in{\DD}^n$
 for $F$ and $G$ of class ${\cal{C}}^1$ and Lipschitz 
from
${\RR}^m$ {\rm[}resp. ${\RR}^n${\rm]} into ${\RR}$, one has $F\circ u\in{\DD}$ and $G\circ v\in{\DD}$ and
$$\textstyle\Gamma[F\circ u,G\circ v]=\sum_{i,j} F_i^{\prime}(u) G_j^{\prime}(v) \Gamma[u_i,v_j]\quad{\PP}{\mbox{-p.s.}}.$$
\indent 3.) the bilinear form ${\cal E}[f,g]={\EE}\Gamma[f,g]$ is closed, i.e. ${\DD}$ is complete under the norm
$\|\,.\,\|_{{\DD}}=(\|\,.\,\|_{L^2({\PP})}^2 +{\cal{E}}[\,.\,,\,.\,])^{\frac{1}{2}}$.

4.)  $1\in{\DD}$ and $\Gamma[1,1]=0$.}

\noindent We always write ${\cal E}[f]$ for ${\cal E}[f,f]$ and $\Gamma[f]$ for $\Gamma[f,f]$.\\

With this definition,  the form ${\cal E}$ defined at point 3.) is a {\it Dirichlet form}. This notion has been
 introduced by A. Beurling and J. Deny as a tool in potential theory, see Beurling and Deny (1958-59), and also Fukushima, Oshima and Takeda (1994). The operator $\Gamma$
is the carr\'{e} du champ or squared field operator associated with ${\cal E}$, it has been studied by several authors in more general
 contexts, see Dellacherie and Meyer (1987), Bouleau and Hirsch (1991). Here we refer to $\Gamma$ as the {\it quadratic error operator} of the error structure. Its intuitive meaning is the
conditional variance of the error.

\noindent{\sc  Example 2.1.}  A simple example of error structure is the term
 $$ ({\RR},{\cal B}({\RR}), m, H^1(m),\gamma)$$ where
$m$ is the normal  law $N(0,1)$ and $$H^1(m)=\{f\in L^2 (m) : f^\prime {\mbox{ in the sense of distributions }}\in L^2(m)\}$$
with $\gamma[f]=f^{\prime 2}$ for $f\in H^1(m)$. This structure is associated to the real valued Ornstein-Uhlenbeck process.
It models an erroneous quantity, say $X$, with normal law whose error does not depends on the value of $X$. Instead, the operator $\gamma[f](x)=f^{\prime 2}(x)x^2$
would model an error proportional to $X$, that is, in the sense of physicists, a constant proportional error. The intuitive relation giving the interpretation of
the quadratic error operator $\gamma$ is 
$$\gamma[f](x)=\EE[(\mbox{error on }f(X))^2|X\!=\!x].$$

\noindent{\thesubsection. How proceeds an error calculation}\setcounter{subsection}{5}\\

\label{subsec:2.4}
Let us suppose we are drawing a triangle with a graduated rule and a protractor: we
 take the polar angle of $OA$
say $\theta_1$, we put $OA=\ell_1$, then we take the angle $(OA,AB)$ say $\theta_2$, and we put $AB=\ell_2$.

1) {\it Choose hypotheses on errors}

$\ell_1,\;\ell_2$ and $\theta_1,\;\theta_2$ and their 
errors can be modeled by the following probability space and operator : :
$$\textstyle((0,L)^2\times(0,\pi)^2, {\cal B}((0,L)^2\times(0,\pi)^2),\frac{d\ell_1}{L}\frac{d\ell_2}{L}\frac{d\theta_1}{\pi}\frac{d\theta_2}{\pi},{\DD},\Gamma)$$
where ${\DD}=\{f\in L^2(\frac{d\ell_1}{L}\frac{d\ell_2}{L}\frac{d\theta_1}{\pi}\frac{d\theta_2}{\pi}) :
\frac{\partial f}{\partial \ell_1}, \frac{\partial f}{\partial \ell_2}, \frac{\partial f}{\partial \theta_1}, 
\frac{\partial f}{\partial \theta_2} \in L^2(\frac{d\ell_1}{L}\frac{d\ell_2}{L}\frac{d\theta_1}{\pi}\frac{d\theta_2}{\pi})\}$ and 
$$\textstyle\Gamma[f]=\ell_1^2(\frac{\partial f}{\partial \ell_1})^2+\ell_1\ell_2\frac{\partial f}{\partial \ell_1}\frac{\partial f}{\partial \ell_2}
+\ell_2^2(\frac{\partial f}{\partial \ell_2})^2+
(\frac{\partial f}{\partial \theta_1})^2+\frac{\partial f}{\partial \theta_1}\frac{\partial f}{\partial \theta_2}
+(\frac{\partial f}{\partial \theta_2})^2,$$ It is easily checked that assumptions 1) 2) 3) 4) of definition 2.1 are fulfilled.

2) {\it Compute the errors on significant quantities using the functional calculus on $\Gamma$}

 \noindent For the coordinates of the point  $B$ for example we have : 
$$\begin{array}{rl}
X_B=&\ell_1\cos \theta_1 +\ell_2\cos(\theta_1+\theta_2)\quad Y_B=\ell_1\sin \theta_1 +\ell_2\sin(\theta_1+\theta_2)\\
\Gamma[X_B]=&\ell_1^2+\ell_1\ell_2(\cos\theta_2+2\sin\theta_1\sin(\theta_1\!+\!\theta_2))+\ell_2^2(1+2\sin^2(\theta_1\! 
+\!\theta_2))\\
\Gamma[Y_B]=&\ell_1^2+\ell_1\ell_2(\cos\theta_2+2\cos\theta_1\cos(\theta_1\!+\!\theta_2))+\ell_2^2(1+2\cos^2(\theta_1\! 
+\!\theta_2))\\
\Gamma[X_B,Y_B]=&-\ell_1\ell_2\sin(2\theta_1\!+\!\theta_2)-\ell_2^2\sin(2\theta_1\!+\!2\theta_2)\\
\end{array}$$ 
Then, according to the problem,  we can for example compute the covariance
of the errors on the area and on the perimeter of the triangle, etc., or obtain that the proportional error on the area  $\Gamma[\mbox{area}]/\mbox{area}^2$ is maximal
for $\theta_2=\pi/2$, etc.

\noindent {\sc Remark 2.1}. If we limit our investigation to variances of the errors, that is to computation with $\Gamma$, then the choice of the {\it a priori} laws is not so crucial
as it could be thought because these computations are done almost surely (using property 2 of definition 2.1. If we are, instead, interested also in biases,
then the {\it a priori} laws are precisely relevant. Biases are represented by an operator which is the generator of the semi-group canonically associated 
with the error structure (see Bouleau and Hirsch (1991) chapter I). In example 2.1 it is $A[u](x)=\frac{1}{2}u^{\prime\prime}(x)-\frac{1}{2}xu^\prime(x)$. If we change the probability
measure $m$ into $f.m$ ($f$ regular), $\gamma$ being unchanged, then the operator $A$ becomes $\tilde{A}[u]=A[u]+\frac{1}{2f}\gamma[f,u]$. 
The operator $u\rightarrow\frac{1}{2f}\gamma[f,u]$ is first order. Absolutely continuous changes of the probability measure $m$ correspond to changes of 
the drift of the bias operator $A$, a variant of Girsanov theorem.\\

\noindent{\thesubsection. Comparison of approaches}\setcounter{subsection}{6}\\

\label{subsec:2.5}
Before looking at the infinite dimensional examples needed in finance, let us try to give an outlook over the different
approaches to error calculus.

At the extreme right-hand side of the table we have the usual probability calculus in which the errors are random variables. 
The knowledge of the joint laws of the quantities and their errors is supposed to be yielded by statistical methods. The errors are finite, 
the propagation of the errors needs computation of image probability laws.

\begin{table}[h]
\label{tab:1}   
\noindent\begin{tabular}{||p{3cm}||p{3cm}|p{3cm}||p{3cm}||} \hline
\multicolumn{1}{||c||}{deterministic}  & \multicolumn{3}{c||}{probabilistic  approaches}\\
 \multicolumn{1}{||c||}{approach} & \multicolumn{3}{c||}{}\\ \hline
\multicolumn{1}{||c||}{Sensitivity}  &\multicolumn{2}{c||}{Extended error calculus using }&  \multicolumn{1}{|c||}{Probability} \\ 
\multicolumn{1}{||c||}{calculus:}&\multicolumn{2}{c||}{Dirichlet forms}&  \\ 
derivation with respect to the parameters of the model & first order calculus only dealing with variances &
second order calculus with variances and biases & \multicolumn{1}{|c||}{theory}\\ \hline
\multicolumn{3}{||c||}{infinitesimal errors} &  \multicolumn{1}{|c||}{ finite errors }\\ \hline
\end{tabular}
\caption{Main classes of error calculi}
\end{table}

At the extreme left-hand side the usual sensitivity calculus consists of computing derivatives with respect to parameters. Let us remark that it applies also
to functional coefficients using Fr\'echet or G\^{a}teaux derivatives.

Between these two purely probabilistic and purely deterministic approaches lies the extended error calculus based on Dirichlet forms.
It supposes the errors infinitely small but takes in account some features of the probabilistic approach allowing to put the computations
and the arguments inside a powerful mathematical theory: the theory of Dirichlet forms. In the same framework
can be performed either a first order calculus on variances which is simple and significant enough for most applications
or a second order calculus dealing with both variances and biases.\\

\noindent{\thesubsection. Main features of the method}\setcounter{subsection}{7}\\

\label{subsec:2.6}
As above in the finite dimensional case of the  triangle, the construction of an error structure on an infinite dimensional stochastic model
is done in two steps
	
1) If there are, as usually, deterministic parameters which can be erroneous or with respect to which a sensitivity
 is wished, these parameters have to be randomized with {\it a priori} laws.

2) Errors operators must be chosen to act on random quantities (initially random or randomized parameters) in order to describe errors, in such
a way that we obtain mathematically an {\it error structure}. 

Several properties of error structures make it easier such a construction.

1) The operation of {\it taking the image} of an error structure by a mapping is quite natural and gives an error structure as soon as the
mapping, even non injective, satisfies some rather weak conditions. In particular if $(\Omega, {\cal A}, {\PP}, {\DD}, \Gamma)$ is an error
structure and if $X$ is a random variable with values in ${\RR}^d$ whose components are
 in ${\DD}$, $({\RR}^d, {\cal B}({\RR}^d),{\PP}_X,{\DD}_X,\Gamma_X)$
is an error structure where ${\PP}_X$ is the law of $X$ and
$$\begin{array}{rcl}
{\DD}_X&=&\{f\in L^2({\PP}_X) : f\circ X\in{\DD}\}\\
\Gamma_X[f]&=&{\EE}[\Gamma[f\circ X]|X\!=\!x],\quad f\in{\DD}.\end{array}$$
\indent 2) If $f\in{\DD}$ and $F$ is Lipschitz from ${\RR}$ to ${\RR}$ then 
$F\circ f\in{\DD}$ and $\Gamma[F\circ f]\leq k\Gamma[f]$. 
For example the structure of example 1 $ ({\RR},{\cal B}({\RR}), m, H^1(m),\gamma)$
possesses an image by the map $x\rightarrow |\sin\sqrt{1+|x|}|$ which is an error structure on $[0,1]$. Such a use of non injective
functions is tricky in the deterministic sensitivity calculus.
More generally if $F$ is a contraction from ${\RR}^d$ into ${\RR}$
 in the following sense
$$|F(x)-F(y)|\leq {\textstyle\sum}_{i=1}^d|x_i-y_i|$$
then for $f_1,f_2,\ldots,f_d\in{\DD}$ one has $F(f_1,f_2,\ldots,f_d)\in{\DD}$ and 
$$\Gamma[F(f_1,f_2,\ldots,f_d)]^\frac{1}{2}\leq{\textstyle\sum}_{i=1}^d\Gamma[f_i]^\frac{1}{2}.$$
This property allows to consider more general images with values in metric spaces as soon as a
 suitable density property is preserved, see Bouleau-Hirsch (1991) chapter V \S 1.3 p 197.

3) The {\it product} of two or countably many error structures is an error structure. It is the mathematical expression
of the independence of the random variables and the non-correlation of the errors. By this way error structures on infinite dimensional
spaces are easily obtained, e.g. on the Wiener space, as we will see in the next part, or on the general Poisson space or other
spaces of stochastic processes, see  Bouleau and Hirsch (1991), Ma and Roeckner (1992), Bouleau (1995).

For later reference we give the following statement.

\noindent{\sc Theorem 2.1}{\it . Product structures

Let $S_n=(\Omega_n,{\cal F}_n, m_n,{\DD}_n, \Gamma_n)$, $n\geq 0$ be error structures.

The term  $S=(\Omega,{\cal F}, m,{\DD},\Gamma)$ defined below is an error structure denoted
 $S=\prod_{n=1}^\infty S_n$ and called the {\rm product structure} of the $S_n$:
$${\textstyle(\Omega,{\cal F}, m)=(\prod_{n=0}^\infty \Omega_n,\bigotimes_{n=0}^\infty {\cal F}_n,\prod_{n=0}^\infty m_n)}$$
$$
\begin{array}{rl}
{\DD}=\{f\in L^2(m) : & \forall n, {\mbox{ for }m\mbox{-a.e. }} \omega=(\omega_0,\omega_1,\ldots) \\
&\mbox{the function }x\rightarrow f(\omega_0,\ldots,\omega_{n-1},x,\omega_{n+1},\ldots)\in{\DD}_n\\
&\mbox{and } \int \sum_n \Gamma_n[f]\;dm<+\infty\}
\end{array}$$ 
and for $f\in{\DD}\quad\Gamma[f]=\sum_n \Gamma_n[f].$
}

Thanks to these properties, is possible the construction of a variety of error structures
 on a given probabilistic model. Now for a rational treatment of a practical case these error hypotheses should be 
obtained by statistical methods. This is connected with the {\it Fisher information theory}, see Bouleau (2001). Anyhow, these statistical
methods are not yet sufficiently studied to be exposed here, especially in the infinite dimensional case we have to use in finance.
Thus we limit ourselves to error computations with {\it a priori} errors chosen the most likely we can. We 
will see that it is significant already.

\setcounter{section}{3}
\begin{center}
{\thesection. ERROR STRUCTURES ON THE WIENER SPACE}
\end{center}
\label{sec:3}
\setcounter{equation}{0}
\setcounter{subsection}{1}
Let us first recall the classical construction of the Brownian motion using the Wiener integral.\\

\noindent{\thesubsection. The Wiener space as Gaussian product space}\setcounter{subsection}{2}\\

\label{subsec:3.1}
Since we aim here at applications we will consider only the case where a  measured space
$(E,{\cal E}, \mu)$ is given which is either $({\RR}_+, {\cal B}({\RR}_+),dt)$ or $([0,1],{\cal B}([0,1]),dt)$ and
 a one-dimensional Brownian motion (for the abstract
 Wiener space setting see Bouleau and Hirsch (1991)).

Let $(\chi_n)$ be an orthonormal basis of $L^2(E,{\cal E}, \mu)$ and let $(g_n)$ be a sequence of i.i.d. reduced Gaussian variables defined on
a probability space $(\Omega, {\cal A}, {\PP})$. To each $f\in L^2(E,{\cal E}, \mu)$ we associate
 $I(f)\in L^2(\Omega, {\cal A}, {\PP})$
by
$$I(f)={\textstyle\sum}_n<f,\chi_n>g_n.$$
then $I$ is an isometric homomorphism from the Hilbert space $L^2(E,{\cal E}, \mu)$ into the Hilbert space
 $L^2(\Omega, {\cal A}, {\PP})$.
If $f$ and $g$ are orthogonal in $L^2(E,{\cal E}, \mu)$, $I(f)$ and $I(g)$ are independent Gaussian random variables and putting
\begin{equation}
B_t={\textstyle\sum}_n <1_{[0,t]}, \chi_n>g_n\quad\quad(t\in[0,1]{\mbox{ or }}t\in{\RR}_+)
\end{equation}
defines a Gaussian stochastic process which is easily shown to be a standard Brownian motion. 
By extending the case where $f$ is a step
function, the random variable $I(f)$ is denoted by ${\textstyle\int} f(s)\;dB_s$
and defines the Wiener integral of $f$.
In this construction we can suppose the space $(\Omega, {\cal A}, {\PP})$ be a product space:
$$(\Omega, {\cal A}, {\PP})=({\RR},{\cal B}({\RR}),m)^{\NN}\quad\quad m=N(0,1)$$
and the $g_n$'s be the coordinate mappings. Thus $\omega=(\omega_0,\ldots, \omega_n,\ldots)$ and $g_n(\omega)=\omega_n$. 

By the functional calculus, as soon as ${\DD}$ and $ \Gamma$ define an error structure on  $(\Omega, {\cal A}, {\PP})$,
say $(\Omega, {\cal A}, {\PP},{\DD},\Gamma)$ for which the $g_n$'s are in ${\DD}$,
this structure is determined by the quantities
\begin{equation}
\Gamma[{\textstyle\int} f(s)\;dB_s]\quad\quad f\in D_0 \mbox{ dense in } L^2[0,1]\mbox{ (resp. } L^2[{\RR}_+])
\end{equation}
because
 it follows that if $F\in{\cal C}^1\cap Lip({\RR}^k)$
$$
\begin{array}{rl}
{\textstyle
\Gamma[F(\int f_1dB, \ldots,\int f_kdB)]}=&\\
{\textstyle\sum_{i,j=1}^k F'_i(\int f_1dB, \ldots)}&{\textstyle F'_j(\int f_1dB, \ldots)\Gamma[\int f_idB,\int f_jdB]}
\end{array}$$ and the random variables 
${\textstyle F(\int f_1dB, \ldots,\int f_kdB)}$ for $F\in{\cal C}^1\cap Lip({\RR}^k)$ and $f_i\in D_0$ are a dense
subspace of $L^2({\PP})$.\\

\noindent{\thesubsection. The Orn\-stein-Uhlen\-beck structure}\setcounter{subsection}{3}\\

Taking $\Gamma[\int f(s)\;dB_s]=\|f\|^2_{L^2}$ and $D_0=L^2[0,1]$ gives a {\it closable} structure
which is of the form
$$\textstyle(\Omega, {\cal A}, {\PP},{\DD},\Gamma)=\prod_{n=0}^\infty(({\RR},{\cal B}({\RR}),m,{\bf d}_n,\gamma_n)$$
where each factor
$({\RR},{\cal B}({\RR}),m,{\bf d}_n,\gamma_n)$ is here a copy of $({\RR},{\cal B}({\RR}), m, H^1(m),\gamma)$  the Orn\-stein-Uhlen\-beck
structure of example 2.1. This error structure is induced by the following perturbation of 
the Brownian path :
\begin{equation}
\omega(t)\longrightarrow \sqrt{e^{-\theta}}\,\omega(t)+\sqrt{1-e^{-\theta}}\,\hat{\omega}(t)
\end{equation}
where $\hat{\omega}$ is an independent Brownian motion and $\theta$ a vanishing parameter.\\

\noindent{\thesubsection. Finer structures of product type}\setcounter{subsection}{4}\\

Taking $\Gamma[\int f(s)\;dB_s]=\|f'\|^2_{L^2}$ and $D_0={\cal C}^1[0,1]$ 
or more generally for $q\in{\NN}^\ast$ (or even $q\in{\RR}_+^\ast$) $\Gamma[\int f(s)\;dB_s]=\|f^{(q)}\|^2_{L^2}$ and $D_0={\cal C}^q[0,1]$ 
give also  {\it closable} structures.
For a suitable choice of the basis $(\chi_n)$ they are also of the form
$$\textstyle(\Omega, {\cal A}, {\PP},{\DD},\Gamma)=\prod_{n=0}^\infty
(({\RR},{\cal B}({\RR}),m,{\bf d}_n,\gamma_n).$$

\noindent{\thesubsection. Error structures of generalised Mehler  type}\setcounter{subsection}{5}\\

Let $p_t$ be a strongly continuous symmetric contraction semi-group
on $ L^2({\RR}_+,dt)$ with generator $(a,{\cal D}a)$ ($a$ being a non necessarily local operator), let us consider the associated closed 
positive quadratic form $(\varepsilon,{\cal D}(\sqrt{-a}))$  defined by $\varepsilon[f]=\|\sqrt{-a}f\|^2$
 (a non necessarily Dirichlet form), then the structure on the Wiener space
induced by the formula
$${\textstyle\Gamma[\int_0^\infty f(s)\;dB_s]=\varepsilon[f],\quad\quad f\in {\cal D}(\sqrt{(-a)}})$$
is {\it closable} and thus defines an error structure.

It corresponds to the semigroup $P_t$ on $L^2(\Omega, {\cal A}, {\PP})$ given by
$$
P_tF=\tilde{{\EE}}\left[F\left(\int_0^\infty(p_{t/2}1_{[0,.]})(u)\;dB_u
+\int_0^\infty(\sqrt{I-p_{t}}1_{[0,.]})(v)\;d\tilde{B}_v\right)\right]
$$
where $\sqrt{I-p_{t}}$ is the positive square root of the positive operator $I-p_{t}$ on $ L^2({\RR}_+,dt)$ and $\tilde{B}$ is an
 auxiliary independent Brownian motion.

\noindent{\sc Remark 3.1}. When a financial model is studied by means of a development in series with respect to a small random change in the coefficients 
(like small noise expansion of a stochastic volatility) it is possible to induce from the perturbation an error structure which manages the variances and
 the biases at the limit when the perturbation is infinitely small. We cannot describe the details here of this standard method. That yields often error structures
outside the class of generalized Mehler type, but still defined by the quantities (3.2).\\

\noindent{\thesubsection. The gradient operator and the derivative}\setcounter{subsection}{6}\\

\label{subsec:3.5}
In  any error structure whose space ${\DD}$ is separable,
 we can define
 a {\it gradient} operator $D$ on ${\DD}$ with values in $L^2({\PP},H)$ where $H$ is an auxiliary Hilbert space:
 $D$ is a continuous application from ${\DD}$ into $L^2({\PP},H)$ such that

1) $\forall U,V\in{\DD}\quad<DU,DV>_H=\Gamma[U,V]$

2) $\forall F\in{\cal C}^1\cap Lip({\RR}^d)$, $\forall X\in{\DD}\quad\quad
\textstyle D(F\circ X)=\sum_{i=1}^d F'_i\circ X.DX_i\quad{\PP}{\mbox{-a.s}}.$
Then
 if $U,V\in{\DD}\cap L^\infty$ which is an algebra, it holds 
$D(UV)=DU.V+U.DV$.

\noindent For example in the case of the Orn\-stein-Uhlen\-beck structure, taking $H=L^2(\RR_+)$ gives

. $\forall h\in L^2({\RR}_+)\quad\quad D_{\!ou}[\int h(s)\,dB_s]=h$

. with suitable hypotheses on the adapted processe $H_t$ (see Nualart (1995))

\quad\quad$D_{\!ou}[\int H_s\,dB_s](t)=H_t+\int(D_{\!ou}H_s)(t)\,dB_s$

. If $U\in{\DD}_{\!ou}$ the Clark formula 

\quad\quad$U={\EE}U+\int_0^\infty{\EE}[D_{\!ou}U(t)|{\cal F}_t]\;dB_t.$

\noindent Now a slight variant of the gradient operator,
the notion of `{\it derivative}', is useful when computing 
errors on solutions of stochastic differential equations thanks to the tool of Ito's formula (this notion has been used and studied
by  Feyel and la Pradelle (1989)).

\noindent{\sc Definition 3.1}. {\it
  Let $(\hat{B}_t)_{t\geq 0}$ be an auxiliary independent Brownian motion. For $U\in{\DD}$ the derivative $U^{\#}$ is a random variable
depending on $\omega$ and $\hat{\omega}$ defined by
$$\textstyle U^{\#}=\int_0^\infty(D_{\!ou}U)(\omega,t)\;d\hat{B}_t.$$
}From the properties of the gradient one gets

. $\Gamma[U]=\hat{{\EE}}[U^{\#2}]$

. For $F\in{\cal C}^1\cap Lip\quad\quad(F\circ U)^{\#}=F'\circ U\;.\;U^{\#}$

\quad\quad$(\int H_s\,dB_s)^{\#}=\int H_s\,d\hat{B}_s+\int H^{\#}_s\,dB_s.$\\

\noindent{\thesubsection. The weighted Orn\-stein-Uhlen\-beck case} \setcounter{subsection}{7}\\

Its meaning for financial models is to consider non necessarly time translation invariant perturbations of the underlying stock price.
It is a special case of the generalised Mehler type:
$$\textstyle\Gamma[\int f(t) dB_t]=\int\alpha(t)f^2(t)dt\quad\quad f\in{\cal D}({\RR}_+)\quad\quad \alpha{\mbox{ measurable }}\geq 0$$

. $\forall h\in L^2({\RR}_+, (1+\alpha)dt)\quad\quad D[\int h(s)\,dB_s]=\sqrt{\alpha(t)}h$

. $(\int h(s)\,dB_s)^{\#}=\int\sqrt{\alpha(t)}\,h(t)\,d\hat{B}_t$

. with suitable hypotheses on the adapted processe $H_t$

\quad\quad$D[\int H_s\,dB_s](t)=\sqrt{\alpha(t)}H_t+\int(DH_s)(t)\,dB_s$

\quad\quad$(\int H_s\,dB_s)^{\#}=\int H_s\sqrt{\alpha(s)}\,d\hat{B}_s+\int H^{\#}_s\,dB_s$

. If $U\in{\DD}\cap{\DD}_{\!ou}$ \quad\quad$DU=\sqrt{\alpha}D_{\!ou}U$

\noindent The generator $(A,{\cal D}A)$ of this structure can easily be seen to verify 
$$\textstyle A[\int f\;dB]=-\int\alpha(s)f(s)dB_s\quad\quad f\in{\cal D}({\RR}_+)$$ which permits (see formula (\ref{generateur}) in the concluding remarks) to compute
$A$ on a dense part of ${\cal D}A$.  

In the sequel, we focuse on the Orn\-stein-Uhlen\-beck case, but the two following lemmas are also valid in the weighted Orn\-stein-Uhlen\-beck case, and part III
and IV extend to that case with only minor changes.

\noindent{\sc Lemma 3.1}.{\it
The conditional expectation operators ${\EE}[\,.\,|{\cal F}_t]$ are orthogonal projectors in ${\DD}$ on errors sub-structures
(closed sub-vector-spaces of ${\DD}$ stable by Lipschitz functions).
}

\noindent{\sc Lemma 3.2}.{\it
Under the same hypotheses, let $\Gamma_t$ be defined from $\Gamma$ by
$$\Gamma_t[{\textstyle(\int f(s)\;dB_s)}]=\Gamma[{\textstyle(\int 1_{[0,t]}f(s)\;dB_s)}]$$
and let $U\rightarrow U^{\#_t}$ the derivation operator associated with $\Gamma_t$, then for }$U\in{\DD}$:
$$({\EE}[U|{\cal F}_t])^{\#}={\EE}[U^{\#_t}|{\cal F}_t].$$

\setcounter{section}{4}
\begin{center}
{\thesection. APPLICATION TO FINANCIAL MODELS}
\end{center}
\label{sec:4}
\setcounter{subsection}{1}
\setcounter{equation}{0}

\noindent{\thesubsection. The Black-Scholes case}\setcounter{subsection}{2}\\

\noindent{\it Notation}

\noindent The interest rate for the bond is constant,
the asset $(S_t)_{t\geq 0}$ is modeled as the solution of the equation $dS_t=S_t(\mu dt+\sigma dB_t)$.
For a European option of the form $f(S_T)$, $T$ fixed deterministic time (see Lamberton and Lapeyre (1997)), the value at time $t\in[0,T]$ of the option is
 $V_t=F(t,S_t,\sigma,r)$ with
\begin{equation}
\label{representationintegrale}
 \quad F(t,x,\sigma,r)=e^{-r(T-t)}\int_{{\RR}}f(xe^{(r-\frac{\sigma^2}{2})(T-t)+\sigma
y\sqrt{T-t}})\frac{e^{-\frac{y^2}{2}}}{\sqrt{2\pi}}dy.
\end{equation}
 If $f$ is Borel with linear growth,
the function $F$ is ${\cal
C}^1$ in $t\in[0,T[$, ${\cal C}^2$ and Lipschitz in $x\in]0,\infty[$,
let us put $$\textstyle{\mbox{delta}}_t=\frac{\partial F}{\partial x}(t,S_t,\sigma,r)\quad\quad\quad
{\mbox{gamma}}_t=\frac{\partial^2 F}{\partial x^2}(t,S_t,\sigma,r)$$

 \noindent$F$ satisfies the equation
$
\frac{\partial F}{\partial t}+\frac{\sigma^2 x^2}{2}\frac{\partial^2
F}{\partial x^2}+rx\frac{\partial F}{\partial x}-rF=0.
$\\

\noindent{\it Hypotheses}

Our choice is governed by an aim of simplicity.

 a) The error on $(B_t)_{t\geq 0}$ is represented by the  Orn\-stein-Uhlen\-beck
error structure.

 b) The errors on the initial value $S_0$, on the volatility $\sigma$, on the rate $r$  are `constant proportional errors' in the sense of physicists
:
$$\begin{array}{rl}
\Gamma[\phi(S_0)]=&\phi'^2(S_0)\,S_0^2\\
\Gamma[\psi(\sigma)]=&\psi'^2(\sigma)\,\sigma^2\\
\Gamma[\xi(r)]=&\xi'^2(r)\,r^2
\end{array}
$$
\indent c) We chose {\it a priori} laws  : lognormal laws on $S_0$ and  $\sigma$,  an exponential law on $r$.

d) We suppose $(B_t)_{t\geq 0}$ and the randomized quantities are independent and their errors uncorrelated. (In a more complete study, these independence
and uncorrelation assumptions would have to be relaxed, in particular to express links between errors on the asset $(S_t)$ and on the volatility $\sigma$).

In other words, the error on a regular function 
$F((B_t)_{t\geq 0}, S_0,\sigma, r)$ will be represented by the product error structure i.e.
$$\Gamma[F((B_t)_{t\geq 0}, S_0,\sigma, r)]=
\Gamma_{ou}[F(., S_0,\sigma, r)]+F'^2_{S_0} S_0^2 +F'^2_\sigma\sigma^2+F'^2_rr^2$$
where $\Gamma_{ou}$ is the  Orn\-stein-Uhlen\-beck quadratic error operator.

Actually, the theory tells us that hedging and pricing formulas do not involve the drift coefficient $\mu$. So we may take $\mu=r$, i.e. we work
under the probability ${\PP}$ such that $\tilde{S_t}=e^{-rt}S_t$, the discounted stock price, is a martingale.
Since $S_t=S_0e^{\sigma B_t+(r-\frac{\sigma^2}{2})t}$ we have 
$$\textstyle\Gamma[S_t]=S^2_t\{\sigma^2\int_0^t\alpha(s)ds +(B_t-\sigma t)^2\sigma^2+t^2\}.$$

\noindent{\thesubsection. Errors on the value and the hedge of a European option}\setcounter{subsection}{3}\\

\label{subsec:4.2}
Let us consider an option of the form $f(S_T)$ where $f$ is Lipschitz.

By the independence hypothesis, the errors on $B$, $S_0$, $\sigma$, $r$ can be managed separately. Let us denote $\Gamma_B$,
$\Gamma_0$, $\Gamma_\sigma$, $\Gamma_r$ the corresponding quadratic operators.\\

\noindent a) {\it Error on the value of the option}

The value of the option is $V_t=F(t,S_t,\sigma,r)$ with $F$ given by (\ref{representationintegrale})

 \noindent a1){\it  Error due to $B$}

$B$ being present only in $S_t$, we have $\Gamma_B[V_t]=(\frac{\partial F}{\partial x}(t,S_t,\sigma,r))^2\Gamma_B[S_t]$ so
\begin{equation}
\begin{array}{rl}
\Gamma_B[V_t]=&{\mbox{delta}_t}^2\;\Gamma_B[S_t]\\
\Gamma_B[V_s,V_t]=&{\mbox{delta}_s}{\mbox{delta}_t}\;\Gamma_B[S_s,S_t]
\end{array}
\end{equation}
with $\Gamma_B[S_s,S_t]=S_sS_t\sigma^2\,(s\!\wedge\! t)$.

\noindent{\sc Proposition 4.1}.
{\it If $f$ is Lipschitz, $V_t$ is in ${\DD}_B$ and when ${t\uparrow T}$
$$V_t=F(t,S_t,\sigma, r)\rightarrow f(S_T)\quad{\mbox{ in }}{\DD}_B{\mbox{ and }}{\PP}-a.s.$$
$$\Gamma_B[V_t]=(\mbox{\rm delta}_t)^2\Gamma_B[S_t]\rightarrow f'^2(S_T)\Gamma_B[S_T]
\quad{\mbox{ in }}L^1{\mbox{ and }}{\PP}-a.s.$$
}

\noindent{\it proof.} Let us suppose first $f\in{\cal C}^1\cap Lip$. By the relation
$$V_t={\EE}[e^{-r(T-t)}f(S_T)|{\cal F}_t]$$ it follows that $V_t\rightarrow f(S_T)$ in $L^p\quad 1\leq p<\infty$ and a.s.

A computation that we shall do in a more general framework later, and that we do not repeat here, gives
$${V_t}^{\#}=e^{-r(T-t)}{\EE}[f'(S_T)S_T|{\cal F}_t]\sigma \hat{B}_t$$
 thus 
$${V_t}^{\#}\rightarrow f'(S_T)S_T\sigma \hat{B}_T{\mbox{ in }} L^2({\PP},L^2(\hat{\Omega},\hat{{\PP}}))$$
and thanks $f(S_T)^{\#}=f'(S_T)S_T\sigma \hat{B}_T$ we obtain
$$V_t\rightarrow f(S_T)\quad{\mbox{ in }}{\DD}_B{\mbox{ and }}{\PP}-a.s.$$
and
$$\Gamma_B[V_t]=e^{-2r(T-t)}({\EE}[f'(S_T)S_T|{\cal F}_t])^2\sigma^2t\rightarrow f'^2(S_T)\Gamma_B[S_T]$$
in $L^1$ and ${\PP}$-a.s. \hfill\fbox\\

The case $f$ only Lipschitz comes from a special property of the one-dimentional functional calculus in error structures
 (see Bouleau and Hirsch (1991) chapter III prop. 2.1.5), the preceding argument still remains valid.

 \noindent a2){\it Error due to $\sigma$.}

We suppose here $f\in{\cal C}^1\cap Lip$. As $V_t=F(t,S_t,\sigma,r)$ 

$$\textstyle\Gamma_\sigma[V_t]=\{F'_x((t,S_t,\sigma,r)\frac{\partial S_t}{\partial\sigma}+F'_\sigma((t,S_t,\sigma,r)\}^2\;\sigma^2$$
and the computation can be done using the integral representation (\ref{representationintegrale}), puting 
$$\tilde{F}(t,x,\sigma)=e^{r(T-t)}F(t,xe^{-r(T-t)},\sigma,r)$$
and remarking that by  (\ref{representationintegrale}) we have
$$\textstyle
\frac{\partial \tilde{F}}{\partial t}+\frac{\sigma^2x^2}{2}\frac{\partial^2\tilde{F}}{\partial x^2}=0\quad\quad
\frac{\partial \tilde{F}}{\partial \sigma}=-\frac{2(T-t)}{\sigma}\frac{\partial \tilde{F}}{\partial t}\quad\quad
\frac{\partial F}{\partial \sigma}=(T-t)\sigma x^2\frac{\partial^2F}{\partial x^2}
$$
and we obtain
\begin{eqnarray}
\Gamma_\sigma[V_t]&=&\{(T-t)\sigma S_t^2\mbox{gamma}_t+S_t(B_t-\sigma t)\mbox{delta}_t\}^2\,\sigma^2\nonumber\\
&=&\{\mbox{vega}_t+S_t(B_t-\sigma t)\mbox{delta}_t\}^2\,\sigma^2
\end{eqnarray}
One gets immediately, for example, the well-known fact that for two European options of payoffs $f_{(1)}(S_T)$ and $f_{(2)}(S_T)$, an option with payoff 
$a_1f_{(1)}(S_T)+a_2f_{(2)}(S_T)$ would have a value $V_0$ at $t=0$ insensitive to $\sigma$, i.e. $\Gamma_\sigma[V_0]=0$, as soon as 
$a_1\,{\mbox{gamma}}_0^{(1)}+a_2\,{\mbox{gamma}}_0^{(2)}=0$.

\noindent a3){\it Error due to $r$.}

We have similarly
$$\Gamma_r[V_t]=\{F'_x((t,S_t,\sigma,r)\frac{\partial S_t}{\partial r}+F'_r((t,S_t,\sigma,r)\}^2\;r^2$$
thus
\begin{equation}\Gamma_r[V_t]=\{tS_t\mbox{delta}_t+{\mbox{rho}}_t\}^2\;r^2.\end{equation}
As a consequence, given several options of payoffs $f_{(i)}¡(S_T)$, $i=1,\ldots, k$, the option of payoff $\sum_i a_if_{(i)}¡(S_T)$ has a value at $t\!=\!0$ insensitive
 to both $\sigma$ and $r$ (i.e. $\Gamma_\sigma[V_0]=\Gamma_r[V_0]=0$) if the vector $a=(a_i)$ is orthogonal to the two vectors
 $({\mbox{gamma}}_0^{(i)})$ and $({\mbox{rho}}_0^{(i)})$.

b) {\it Error on the hedging portfolio}

Here we limit ourselves to the error due to $(B_t)$. We suppose $f$ and $f'$ in ${\cal C}^1\cap Lip$.
The hedging equation is 
$${\textstyle e^{-rt}F(t,S_t,\sigma,r)=F(0,S_0,\sigma,r)+\int_0^t H_s\;d\tilde{S}_s}$$
where the adapted process $H_t$ is the quantity of stock in the portfolio :
$$H_t=\mbox{delta}_t=\frac{\partial F}{\partial x}(t,S_t,\sigma,r)=e^{-r(T-t)}{\EE}[f'(S_T)S_T|{\cal F}_t]\frac{1}{S_t}.$$
By the same method as for $V_t$ we obtain 
\begin{equation}
\begin{array}{rl}
\Gamma_B[H_t]=&({\mbox{gamma}}_t)^2\Gamma_B[S_t]\\
\Gamma_B[H_s,H_t]=&{\mbox{gamma}}_s{\mbox{gamma}}_t\Gamma_B[S_s,S_t]
\end{array}
\end{equation}
{\sc Proposition 4.2}. {\it
If $f,f'\in{\cal C}^1\cap Lip$, then $H_t\in {\DD}$ and as $t\uparrow T$
$$\begin{array}{rl}
H_t\rightarrow& f'(S_T)\quad \mbox{ in }{\DD}_B\mbox{ and a.s.}\\
\Gamma_B[H_t]\rightarrow& f''^2(S_T)\Gamma_B[S_T]\quad \mbox{ in }L^1({\PP})\mbox{ and a.s.}
\end{array}$$
}

\noindent{\sc Remark 4.1}. These results show that the Greeks introduced by practioners have a direct sense as 
sensitivity of the value $V_t$ and of the hedging $H_t$ to  perturbations. This is of course not surprising, the method makes more
precise the correlations of errors. It gives also a tool to study the absolute continuity of joint laws as we explain now.

The preceding computations show easily that in the Black-Scholes model, if $U_1$ and $U_2$ are two random variables taken among the following
quantities defined at a fixed instant $t$: $S_t,\;V_t(f_1),\;V_t(f_2),\;H_t(f_1),\;H_t(f_2)$, then the matrix $\Gamma[U_i,U_j]$ is singular:
 the errors on these quantities are linked. This comes from the fact that  the law of e.g. the pair $(V_t(f_1),V_t(f_2))$ is carried
 by the $\lambda$-parametrized curve:
$$\textstyle\begin{array}{rl}
y=&\exp{-r(T-t)}P_{T-t}f_1(\lambda)\\
x=&\exp{-r(T-t)}P_{T-t}f_2(\lambda)
\end{array}
$$
where $(P_t)$ is the transition semigroup of $(S_t)$. The same phenomenon  happens in any more general Markovian model.

On the contrary the random quantities involving several different instants have generally non-linked errors. Thus for example
if $U_1=S_T$ and $U_2=\int_0^Te^{-s} H_sS_s ds$ (discounted immobilization of the portfolio) the matrix $\Gamma[U_i,U_j]$ is a.s. regular as soon as 
$f$ is not constant,  hence, by the absolute continuity criterion (Bouleau and Hirsch (1986) or Nualart (1995) thm 2.1.2) the law of the pair 
$(S_T,\;\int_0^T e^{-s}H_sS_s ds)$ possesses a density.\\

\noindent c) {\it More general errors on $(B_t)$}

The relations (4.2) (4.3) (4.4) (4.5) still hold in the weighted Orn\-stein-Uhlen\-beck case, and also, with suitable hypotheses, if we consider 
more general error structures on the Wiener space. Let
 us consider, as mentioned above,  a structure induced by a closed positive quadratic form $\varepsilon$ on $L^2({\RR}_+,dt)$
with 
$${\textstyle\Gamma_B[\int f\;dB]=\varepsilon[f]}$$ for $f$ in the domain of $\varepsilon$ with, for example,
$$\begin{array}{l}
\mbox{i)}\quad\quad\varepsilon[f]=\int_0^\infty \int_0^\infty(f(s)-f(t))^2\beta(s)\beta(t)\;dsdt\\
\mbox{ii)}\quad\quad\varepsilon[f]=\int_0^1(f^{(q)}(s))^2ds
\end{array}$$ where $f^{(q)}$ is the fractional derivative of order $q$,

\noindent then the formulas 
$$\begin{array}{rl}
\Gamma_B[V_t]=&({\mbox{delta}_t})^2\;\Gamma_B[S_t]\\
\Gamma_B[H_t]=&({\mbox{gamma}}_t)^2\Gamma_B[S_t]
\end{array}$$ remain valid as soon as $S_t\in{\DD}$ i.e.

\noindent in case i) if $\beta\in L^1({\RR}_+,dt)$ and $\Gamma_B[S_t]=S_t^2\sigma^22\int_t^\infty
\beta(s)ds\int_\infty^t\beta(s)ds$.

\noindent in case ii) if $q\in(0,\frac{1}{2})$ and $\Gamma_B[S_t]=S_t^2\sigma^2\sum_{n=1}^\infty\frac{4(1-\cos 2\pi nt)}{(2\pi n)^{2(1-q)}}.$

Now in the case
$$\textstyle\mbox{iii)}\quad\quad\varepsilon[f]=\int_0^\infty (\sum_{i=1}^d a_i(s)f^{(i)}(s))^2\;ds$$ we do not 
have anymore $1_{[0,t]}\in\mbox{dom}(\varepsilon)$, hence $B_t$ doesn't belong to ${\DD}$. Such error structures are
more convenient to model errors on processes with finite variation.
\setcounter{section}{5}
\begin{center}
{\thesection. MODELS WITH LEVEL DEPENDENT VOLATILITY}
\end{center}
\setcounter{subsection}{1}
\setcounter{equation}{0}
\label{sec:5}
We will display the method in the case of a complete market, the probability being a martingale measure and
 for a simple one-dimensional diffusion model.

The stock is supposed to be the solution of the equation
$$dX_t=X_t\sigma(t,X_t)\,dB_t+X_tr(t)\,dt.$$
We limit the study to the error due to $(B_t)$ which is  defined by an Orn\-stein-Uhlen\-beck structure:
$${\textstyle\Gamma[\int_0^\infty h(s)\,dB_s]=\int_0^\infty h^2(s)\,ds}.$$ The rate is deterministic, the function $\sigma(t,x)$ will be supposed bounded
with bounded derivative in $x$ uniformly for $t\in[0,T]$.

Let $f(X_T)$ be a European option. Its value at time $t$ is
$${\textstyle V_t={\EE}[\exp(-\int_t^Tr(s)ds)f(X_T)|{\cal F}_t]}$$ the hedging portfolio is given by the adapted
process $H_t$ which satisfies
\begin{equation}
\label{previsibleRepresentation}
{\textstyle\tilde{V}_t=\exp(-\int_0^tr(s)ds)V_t=V_0+\int_0^tH_s\,d\tilde{X}_s}
\end{equation}
where $\tilde{X}_s=\exp(-\int_0^tr(s)ds)X_t$.

We proceed as follows: from the equation
$$X_t=X_0+\int_0^tX_s\sigma(s,X_s)dB_s+\int_0^tr(s)X_sds$$
we obtain
$$X_u^{\#}=\int_0^u(\sigma(s,X_s)+X_s\sigma'_x(s,X_s))X_s^{\#}\,dB_s+\int_0^u
X_s\sigma(s,X_s)\,d\hat{B}_s+\int_0^ur(s)X_s^{\#}ds$$
this equation is solved by putting
$$
\begin{array}{rl}
K_s=&\sigma(s,X_s)+X_s\sigma'_x(s,X_s)\\
M_u=&\exp\left\{\int_0^uK_s\,dB_s-\frac{1}{2}\int_0^uK_s^2ds+\int_0^ur(s)ds\right\}
\end{array}
$$
and remarking that
$$X_u^{\#}=M_u\int_0^u\frac{X_s\sigma(s,X_s)}{M_s}\,d\hat{B}_s.$$
a) Let us first suppose $f\in {\cal C}^1\cap Lip$ and let us define $Y=\exp(-\int_t^Tr(s)ds)f(X_T)$. To compute
$({\EE}[Y|{\cal F}_t])^{\#}$ we apply the second lemma  of section 3:
$${\textstyle Y^{\#_t}=\exp(-\int_t^Tr(s)ds)f'(X_T)X_T^{\#_t}}$$ and 
$${\textstyle({\EE}[Y|{\cal F}_t])^{\#}=\exp(-\int_t^Tr(s)ds){\EE}[f'(X_T)X_T^{\#_t}|{\cal F}_t]}$$
$${\textstyle=\exp(-\int_t^Tr(s)ds)}{\EE}[f'(X_T)M_T|{\cal F}_t]\int_0^t\frac{X_s\sigma(s,X_s)}{M_s}\,d\hat{B}_s$$
and  the second lemma  gives
\begin{equation}\begin{array}{rl}
\Gamma[V_t]=&\Gamma[{\EE}[Y|{\cal F}_t]]\\
=&\exp(-2\int_t^Tr(s)ds)({\EE}[f'(X_T)M_T|{\cal F}_t])^2\int_0^t\frac{X_s^2\sigma^2(s,X_s)}{M_s^2}\,ds
\end{array}
\end{equation}
this yields also the cross error of $V_s$ and $V_t$ which is usefull to compute errors on 
random variables such that $\int_0^Th(s)dV_s$ or $\int_0^Th(s)V_sds$.
\begin{equation}\begin{array}{rl}
\Gamma[V_s,V_t]=&\exp(-\int_s^Tr(s)ds-\int_t^Tr(s)ds)\\
&\quad{\EE}[f'(X_T)M_T|{\cal F}_s]{\EE}[f'(X_T)M_T|{\cal F}_t]\\
&\quad\quad\quad\quad\quad\quad\quad\quad\int_0^{s\wedge t}\frac{X_u^2\sigma^2(u,X_u)}{M_u^2}\,du.
\end{array}
\end{equation}
With our hypotheses as $t\uparrow T$
$$\Gamma[V_t]\rightarrow f'^2(X_T)M_T^2\int_0^T\frac{X_s^2\sigma^2(s,X_s)}{M_s^2}\,ds=f'^2(X_T)\Gamma[X_T]$$
in $L^1({\PP})$ and a.s.

\noindent b) Now to deal with $H_t$, let us remark first that $H_t$ is easily obtained by the Clark formula.
 The formula \ref{previsibleRepresentation}
gives
$$H_t{\textstyle\exp(-\int_0^tr(s)ds)}X_t\sigma(X_t)=D_{ad}[{\textstyle\exp(-\int_0^Tr(s)ds)}f(X_T)]$$
where $D_{ad}$ is the adapted O-U-gradient defined by 
$$D_{ad}[Z](t)={\EE}[DZ(t)|{\cal F}_t].$$
Since
$$D[{\textstyle\exp(-\int_0^Tr(s)ds)}f(X_T)]={\textstyle\exp(-\int_0^Tr(s)ds)}f'(X_T)(DX_T)(t)$$
we have from the computation done for $V_t$
$$D[{\textstyle\exp(-\int_0^Tr(s)ds)}f(X_T)]={\textstyle\exp(-\int_0^Tr(s)ds)}{\EE}[f'(X_T)M_T|{\cal F}_t]\frac{X_t\sigma(t,X_t)}{M_t}.$$
Thus 
$$H_t={\textstyle\exp(-\int_t^Tr(s)ds)}{\EE}[f'(X_T)M_T|{\cal F}_t]\frac{1}{M_t}.$$
Now supposing $f$ and $f'\in{\cal C}^1\cap Lip$ we apply the same method as for obtaining $\Gamma[V_t]$ which leads to
\begin{eqnarray}
\Gamma[H_t]&=&{\textstyle\exp(-2\int_t^Tr(s)ds)}\nonumber\\
&&\quad\left({\EE}[\frac{M_T}{M_t}(f''(X_T)M_T+f'(X_T)Z_t^T|{\cal F}_t]\right)^2\int_0^t\frac{X_u^2\sigma^2(u,X_u)}{M_u^2}du\\
&&\nonumber\\
\mbox{with }\quad Z_t^T&=&\int_t^TL_sdB_s-\int_t^T K_sL_sM_sds\nonumber\\
\mbox{and }\quad K_s&=&\sigma(s,X_s)+X_s\,\sigma^\prime_x(s,X_s)\nonumber\\
L_s&=&2\sigma^\prime_x(s,X_s)+X_s\,\sigma^{\prime\prime}_{x^2}(s,X_s).\nonumber
\end{eqnarray}
If we introduce the following notation which, in our present Markovian model, gives the probabilistic interpretation of the usual Greeks
$$\begin{array}{rl}
\mbox{delta}_t=&H_t={\textstyle\exp(-\int_t^Tr(s)ds)}{\EE}[f'(X_T)M_T|{\cal F}_t]\frac{1}{M_t}\\
\mbox{gamma}_t=&{\textstyle\exp(-\int_t^Tr(s)ds)}{\EE}[\frac{M_T^2}{M_t^2}(f''(X_T)+\frac{M_T}{M_t^2}f'(X_T)Z_t^T|{\cal F}_t]
\end{array}$$
we can summarize some formulas of this case with level dependent volatility by
$$\begin{array}{rlrl}
V_t^{\#}=&\mbox{delta}_t X_t^{\#}&&\\

\Gamma[V_t]=&\mbox{delta}_t^2\Gamma[X_t]&\Gamma[V_s,V_t]=&\mbox{delta}_s\mbox{delta}_t\Gamma[X_s,X_t]\\

&&&\\
H_t^{\#}=&\mbox{gamma}_t X_t^{\#}&\Gamma[H_s,H_t]=&\mbox{gamma}_s\mbox{gamma}_t\Gamma[X_s,X_t]\\
\Gamma[H_t]=&\mbox{gamma}_t^2\Gamma[X_t]&\Gamma[V_s,H_t]=&\mbox{delta}_s\mbox{gamma}_t\Gamma[X_s,X_t]\\
&&&\\
\Gamma[X_t]=&M_t^2\int_0^t\frac{X_u^2\sigma^2(u,X_u)}{M_u^2}du\quad&
\Gamma[X_s,X_t]=&M_sM_t\int_0^{s\wedge t}\frac{X_u^2\sigma^2(u,X_u)}{M_u^2}du.
\end{array}$$
\setcounter{section}{6}
\begin{center}
{\thesection. CONCLUDING REMARKS}
\end{center}
\setcounter{subsection}{1}
\setcounter{equation}{0}
\label{sec:6}
The error calculus based on Dirichlet forms begins at present to be used by modelisators in economics and  finance. It is too early to give an account of its applications.
What we have done is just presentating  how this tool runs through stochastic models including SDE's.

Among the directions of research let us mention that this approach yields {\it new integration by parts formulas} which have been shown to be useful to compute
the Greeks by Monte Carlo methods. Also, it allows to perform Malliavin calculus on the Monte Carlo sample space, that
 is after discretization instead of before. Let us indicate briefly the idea : 

Consider the error structure
$$
\begin{array}{l}
(\Omega,{\cal A},{\PP},{\DD},\Gamma)=((0,1),{\cal B}(0,1), dx,\,\dd\,,\gamma)^{\NN^\ast}\\
\dd = \{u\in L^2(0,1)\quad x(1-x)u^\prime(x)\in H_0^1(0,1)\}\\
\gamma[u](x)=x^2(1-x)^2u^{\prime 2}(x).
\end{array}$$
Let us denote $U_n$ the coordinate maps, taking ${\cal H}=\ell^2$ this structure admits the following gradient : if $F=f(U_1, U_2,\ldots, U_n, \ldots)$
belongs to $\DD$
$$DF=(f^\prime_n(U_1, U_2,\ldots, U_n, \ldots)U_n(1-U_n))_{n\geq 1}$$ and if $a\in \ell^2$ we have the following integration by parts formula
$$\textstyle\EE[<DF,a>_{\ell^2}]=-\EE[F\sum_{n=0}^\infty a_n(1-2U_n)].$$
Let us take for example the following discrete approximation of an SDE: 
$$
S_{n+1}=S_n+\sigma(S_n)(Y_{n+1}-Y_n)\;,\quad\quad
S_0=x,
$$ where $Y_{n+1}-Y_n=\lambda\;\xi(n+1,U_{n+1})$,  we get easily under regularity assumptions on $\xi$ and $\sigma$ :
$$\begin{array}{rl}
\frac{d}{dx}\EE[\Psi(S_N)]=&\EE[\Psi(S_N)(\frac{\xi^{\prime\prime}(1,U_1)(1+\sigma^\prime(x)\xi(1,U_1))}{\sigma(x)\xi^{\prime 2}(1,U_1)}-\frac{\sigma^\prime(x)}{\sigma(x)})]\\
\frac{d}{d\lambda}\EE[\Psi(S_N)]=&-\EE[\Psi(S_N)\sum_{n=1}^N\frac{d}{dU_n}\frac{\xi(n,U_n)}{\xi^\prime(n,U_n)}].
\end{array}$$
\indent As a second direction of research, let us sketch how to do a sensitivity analysis of the solution of an SDE with respect to a functional coefficient.

Let us consider that the level dependent volatility of section 5 is a function in a vector space $\Omega_1$ equipped with an error structure 
$(\Omega_1,{\cal A}_1,{\PP}_1,{\DD}_1,\Gamma_1)$ such that $\PP_1$-almost every function in $\Omega_1$ be of class ${\cal C}^1$ and Lipschitz, and that
the linear form $V_x$ defined by $$V_x(\sigma)=\sigma(x)$$ belong to $\DD_1$. Defining the notation $\sigma^{\#}$ by
$\sigma^{\#}(x)=(V_x)^{\#}(\sigma)$ we can show that if $Y$ is an erroneous independent random variable defined on an other error structure, the following
formula holds
$$(\sigma(Y))^{\#}=\sigma^{\#}(Y)+\sigma^\prime(Y)Y^{\#}.$$ Then similar computations to those of section 5 can be done. For example suppose $\sigma$ is
represented for numerical evaluation as
$$\sigma(t,x)=\sum_{p,q}a_{pq}\,t^px^q$$ where the $a_{pq}$'s are random and erroneous with $\Gamma[a_{pq}]=a_{pq}^2,\;$ $\Gamma[a_{pq},a_{ij}]=0$ if 
$(p,q)=\!\!\!\!\!/\;\,(i,j)$, we obtain that the error on $\sigma$ transfers to $X_t$ in the following way:
$$\Gamma[X_t]=M_t^2\sum_{p,q}\left(\int_0^t\frac{s^pX_s^{q+1}}{M_s}(dB_s-K_s\,ds)\right)^2a_{pq}^2$$ where $M_t,\;K_t$ have the same meaning as in section 5.

 Because estimates of biases are important in financial models especially for pricing, see e.g. Hull and White (1988), let us mention shortly  what would
 be the second order calculus with variances and biases mentionned above in the table 2.1 of
 section 2. In an error
structure $(\Omega,{\cal A},{\PP},{\DD},\Gamma)$ the bias of the error on a random variable $X$ (i.e. the conditional expectation 
of the error) is represented by the generator $A$ of the semi-group canonically associated with the error structure acting on $X$, see Bouleau and Hirsch (1991). 
It has a domain ${\cal D}A$ smaller than ${\DD}$. The functional calculus on $A$ follows the following rules: for
 all $F\in{\cal C}^2({\RR}^d)$, $\forall f_i$ locally in $ {\cal D}A$, $i=1,\ldots,d$, then $F(f_1,\ldots,f_d)$ is locally in ${\cal D}A$ and 
\begin{equation}
\label{generateur}A[F(f)]=\sum_{i=1}^dF'_i(f)Af_i+\frac{1}{2}\sum_{i,j=1}^dF''_{ij}(f)\Gamma[f_i,f_j]
\end{equation}
On the Black-Scholes model, with the  O-U hypotheses and concerning solely the error due to $(B_t)$, we obtain :
$$\begin{array}{l}
\Gamma[B_t]\!=t,\quad\quad\quad  \Gamma[S_t]=S_t^2 \sigma^2 t,  \quad\quad \quad
A[B_t]\!=-B_t,\\
 A[S_t]=-S_t\sigma B_t+\frac{1}{2}\sigma^2 S_t t,\quad\quad \quad
A[V_t]\!=\mbox{delta}_tA[S_t]+\frac{1}{2}\mbox{gamma}_t\Gamma[S_t],\\
A[H_t]\!=\mbox{gamma}_t A[S_t]+\frac{1}{2}\frac{\partial^3 F}{\partial x^3}(t,S_t,\sigma,r)\Gamma[S_t].
\end{array}$$ 
\indent Except at time $t=0$ (since the perturbation 3.3 doesn't move at $t=0$ but we can imagine it starts farther in the past) we see that an error on the path
of the stock $(S_t)_{t\geq 0}$ induces biases on the price $V_t$
and on the hedge $H_t$ (involving the `Greek' $\frac{\partial^3F}{\partial x^3}$). May these biases due to a fuzzy stock price be an interpretation of the bid-ask ?
This interesting question needs certainly  more complete investigations since when a transaction occurs the price is erroneous
while the amount of stock which is bought or sold is not. Anyhow it is an auspicious project to 
 understand and to take in account the consequences of the bid-ask on pricing and hedging procedures of  tractable financial models with the help of
 the error calculus on variances (operator $\Gamma$) and biases (operator $A$).

\begin{list}{}
{\setlength{\itemsep}{0cm}\setlength{\leftmargin}{0.5cm}\setlength{\parsep}{0cm}\setlength{\listparindent}{-0.5cm}}
  \item\begin{center}
{\small REFERENCES}
\end{center}\vspace{0.4cm}

 {\sc Beurling, A., Deny, J.} (1958-59): Espaces de Dirichlet, I. le cas \'{e}l\'{e}\-mentaire, {\it Acta Math.} 99, 203-224 
; Dirichlet spaces, {\it Proc. Nat. Acad.  Sci. U.S.A.} 45, 206-215.
 
{\sc Bogachev, V. I., Roeckner, M.} (1995): Mehler formula and capacities for infinite dimensional Orn\-stein-Uhlen\-beck processes with general
linear drift, {\it Osaka J. Math.} 32, 237-274.

{\sc Bouleau, N.} (1995): Construction
 of Dirichlet structures, {\it in : Potential theory ICPT 1994}, de Gruyter.
 
 {\sc Bouleau, N.} (2001): Calcul d'erreur complet lipschitzien et formes de Dirichlet, {\it J. Math. pures et
 appl.} 80, 9, 961-976.

 {\sc Bouleau, N.,} and {\sc Hirsch, F.} (1986): Propri\'{e}t\'{e}s d'absolue continuit\'{e} dans les espaces de Dirichlet et application
aux EDS, {\it in S\'{e}m. probabilit\'{e} XX, Lecture notes in Mathematics} 1204, Springer, 131-161.

 {\sc Bouleau, N.,} and {\sc Hirsch, F.} (1991): {\it Dirichlet forms and analysis on Wiener space,} De Gruyter.

{\sc Bouleau, N.,} and {\sc Lamberton, D.} (1989): Residual risks and hedging strategies in Markovian markets, {\it Stochast. Process. Appl.}
33, 131-150.

{\sc Dellacherie,} and {\sc Cl., Meyer, P. A.} (1987): {\it Probabilit\'{e}s et potentiel.} Hermann.

 {\sc Feyel, D.,} and {\sc la Pradelle, A. de}  (1989): Espaces de Sobolev Gaussiens, {\it Ann. Inst. Fourier}, 39-4, 875-908.
 
{\sc Fourni\'e, E., Lasry, J.-M., Lebuchoux, J., Lions, P.-L.,} and {\sc Touzi, N.} (1999):
Application of Malliavin calculus to Monte Carlo methods in finance, {\it Finance and Stochastics}, 391-412.

{\sc Fukushima, M., Oshima, Y.,} and {\sc Takeda, M.} (1994): {\it Dirichlet forms and Markov processes.} De Gruyter.

{\sc Hull, J.} and {\sc White, A.} (1988) An analysis of the bias in option pricing caused by a stochastic volatility,
 {\it Advances in Futures and Options Research}, vol 3, 29-61.

 {\sc Lamberton, D.,} and {\sc Lapeyre, B.} (1997): {\it Introduction au calcul stochastic appliqu\'{e} \`{a} la finance}. Ellipses.

{\sc Ma, Z.,} and {\sc Roeckner, M.} ( 1992): {\it Dirichlet forms}. Springer.

 {\sc Malliavin, P.} (1997): {\it Stochastic analysis}. Springer.

 {\sc Nualart, N.} (1995): {\it 
The Malliavin calculus and related topics}. Springer.
\end{list}

\end{document}